\documentclass[12pt]{amsart}

\usepackage{amssymb, amsthm, bm}
\usepackage{fullpage}

\newtheorem{theorem}{Theorem}
\newtheorem{corollary}{Corollary}
\newtheorem{lemma}{Lemma}

\theoremstyle{definition}

\newcommand{\eps}{\varepsilon}
\newcommand{\sfrac}[2]{{\textstyle \frac {#1}{#2}}}

\newcommand{\pmodulo}[1]{\;(\mathrm{mod}\;#1)}

\title{On sums of primes from Beatty sequences}
\author{Angel V. Kumchev}
\thanks{The author's research was partially supported by a Towson University summer research fellowship}
\address{Department of Mathematics, Towson University, Towson, MD 21252-0001, U.S.A.}
\email{akumchev@towson.edu}

\begin{document}

\begin{abstract}
  Let $k \ge 2$ and $\alpha_1, \beta_1, \dots, \alpha_k, \beta_k$ be reals such that the $\alpha_i$'s are irrational and greater than $1$. Suppose further that some ratio $\alpha_i/\alpha_j$ is irrational. We study the representations of an integer $n$ in the form
  \[
    p_1 + p_2 + \dots + p_k = n,
  \]
  where $p_i$ is a prime from the Beatty sequence
  \[
    \mathcal B_i = \left\{ n \in \mathbb N : n = \left[ \alpha_i m + \beta_i \right] \text{ for some } m \in \mathbb Z \right\}.
  \]
\end{abstract}

\maketitle

\section{Introduction}

Ever since the days of Euler and Goldbach, number-theorists have been fascinated by additive representations of the integers as sums of primes. The most famous result in this field is I.M. Vinogradov's three primes theorem \cite{IVin37a}, which states that every sufficiently large odd integer is the sum of three primes. Over the years, a number of authors have studied variants of the three primes theorem with prime numbers restricted to various sequences of arithmetic interest. For instance, a recent work by Banks, G\"ulo\u glu and Nevans \cite{BaGuNe07} studies the question of representing integers as sums of primes from a Beatty sequence. Suppose that $\alpha$ and $\beta$ are real numbers, with $\alpha > 1$ and irrational. The \emph{Beatty sequence} $\mathcal B_{\alpha, \beta}$ is defined by
\[
  \mathcal B_{\alpha, \beta} = \left\{ n \in \mathbb N : n = \left[ \alpha m + \beta \right] \text{ for some } m \in \mathbb Z \right\}.
\]
(Henceforth, $[\theta]$ represents the integer part of the real number $\theta$.) Banks \emph{et al.} proved that if $k \ge 3$, then every sufficiently large integer $n \equiv k \pmodulo 2$ can be expressed as the sum of $k$ primes from the sequence $\mathcal B_{\alpha, \beta}$, provided that $\alpha < k$ and $\alpha$ ``has a finite type'' (see below). In their closing remarks, the authors of \cite{BaGuNe07} note that their method can be used to extend the main results of \cite{BaGuNe07} to representations of an integer $n$ in the form 
\begin{equation}\label{1}
  p_1 + p_2 + \dots + p_k = n,
\end{equation}
where $p_i \in \mathcal B_{\alpha, \beta_i}$. However, they remark that ``for a sequence $\alpha_1, \dots, \alpha_k$ of irrational numbers greater than $1$, it appears to be much more difficult to estimate the number of representations of $n \equiv k \pmodulo 2$ in the form \eqref{1}, where $p_i$ lies in the Beatty sequence $\mathcal B_{\alpha_i, \beta_i}$.'' The main purpose of the present note is to address the latter question in the case when at least one of the ratios $\alpha_i/\alpha_j$, $1 \le i,j \le k$, is irrational.

Let $\alpha_1, \beta_1, \dots, \alpha_k, \beta_k$, $k \ge 2$, be real numbers, and suppose that $\alpha_1, \dots, \alpha_k$ are irrational and greater than $1$. For $i = 1, \dots, k$, we denote by $\mathcal B_i$ the Beatty sequence $\mathcal B_{\alpha_i, \beta_i}$. We write 
\begin{equation}\label{2}
  R(n) = R(n; \bm{\alpha}, \bm{\beta}) = \sum_{\substack{ p_1 + \dots + p_k = n\\ p_i \in \mathcal B_i}} (\log p_1) \cdots (\log p_k),
\end{equation}
where the summation is over the solutions of \eqref{1} in prime numbers $p_1, \dots, p_k$ such that $p_i \in \mathcal B_i$. Similarly to
\cite{BaGuNe07}, we shall use the Hardy--Littlewood circle method to obtain an asymptotic formula for $R(n)$. The circle method requires some quantitative measure of the irrationality of the $\alpha_i$'s in the form of hypotheses on the rational approximations to the $\alpha_i$'s. Let $\| \theta \|$ denote the distance from the real number $\theta$ to the nearest integer. We say that an $s$-tuple $\theta_1, \dots, \theta_s$ of real numbers \emph{is of a finite type}, if there exists a real number $\eta$ such that the inequality
\begin{equation}\label{3}
  \| q_1\theta_1 + \dots + q_s\theta_s \| < \max( 1, |q_1|, \dots, |q_s| )^{-\eta}
\end{equation}
has only finitely many solutions in $q_1, \dots, q_s \in \mathbb Z$. In particular, the reals in an $s$-tuple of a finite type are irrational and linearly independent over $\mathbb Q$. Our main result can now be stated as follows.

\begin{theorem}\label{th1}
  Let $k \ge 3$ and let $\alpha_1, \beta_1, \dots, \alpha_k, \beta_k$ be real numbers, with $\alpha_1, \dots, \alpha_k > 1$. Suppose that each individual $\alpha_i$ is of a finite type and that at least one pair $\alpha_i^{-1}, \alpha_j^{-1}$ is also of a finite type. Then, for any fixed $A > 0$ and any sufficiently large integer $n$, one has
  \begin{equation}\label{4}
    R(n; \bm{\alpha}, \bm{\beta}) = \frac {\mathfrak S_k(n)n^{k - 1}}{\alpha_1 \cdots \alpha_k(k - 1)!} + O \big( n^{k - 1}(\log n)^{-A} \big),
  \end{equation}
  where $R(n; \bm{\alpha}, \bm{\beta})$ is the quantity defined in \eqref{2} and $\mathfrak S_k(n)$ is given by
  \begin{equation}\label{5}
    \mathfrak S_k(n) = \prod_{p \mid n} \left( 1 + \frac {(-1)^k}{(p - 1)^{k - 1}} \right) \prod_{p \nmid n} \left( 1 + \frac {(-1)^{k + 1}}{(p - 1)^k} \right).
  \end{equation}
  The implied constant in \eqref{4} depends at most on $A, \bm{\alpha}, \bm{\beta}$.
\end{theorem}

Since $1 \ll \mathfrak S_k(n) \ll 1$ when $n \equiv k \pmodulo 2$ and $\mathfrak S_k(n) = 0$ when $n \not\equiv k \pmodulo 2$, Theorem~\ref{th1} has the following direct consequence.

\begin{corollary}\label{c1}
  Let $k \ge 3$ and suppose that $\alpha_1, \beta_1, \dots, \alpha_k, \beta_k$ are real numbers subject to the hypotheses of Theorem \ref{th1}. Then, every sufficiently large integer $n \equiv k \pmodulo 2$ can be represented in the form \eqref{1} with $p_i \in \mathcal B_i$, $1 \le i \le k$.
\end{corollary}

After some standard adjustments, the techniques used in the proof of Theorem \ref{th1} yield also the following result on sums of two Beatty primes.

\begin{theorem}\label{th2}
  Suppose that $\alpha_1, \beta_1, \alpha_2, \beta_2$ are real numbers, with $\alpha_1, \alpha_2 > 1$. Suppose further that the pair $\alpha_1^{-1}, \alpha_2^{-1}$ is of a finite type. Then, for any fixed $A > 0$, and all but $O( x(\log x)^{-A} )$ integers $n \le x$, one has
  \[
    R(n; \bm{\alpha}, \bm{\beta}) = (\alpha_1\alpha_2)^{-1}\mathfrak S_2(n)n + O \big( n(\log n)^{-A} \big),
  \]
  where $R(n; \bm{\alpha}, \bm{\beta})$ is the quantity defined in \eqref{2} and $\mathfrak S_2(n)$ is given by \eqref{5} with $k = 2$. The implied constants depend at most on $A, \bm{\alpha}, \bm{\beta}$.
\end{theorem}

Since for even $n$, $1 \ll \mathfrak S_2(n) \ll \log\log n$, we have the following corollary to Theorem \ref{th2}.

\begin{corollary}\label{c2}
  Suppose that $\alpha_1, \beta_1, \alpha_2, \beta_2$ are real numbers subject to the hypotheses of Theorem \ref{th2}. Then, for any fixed $A > 0$, all but $O( x(\log x)^{-A} )$ even integers $n \le x$ can be represented as sums of a prime $p_1 \in \mathcal B_1$ and a prime $p_2 \in \mathcal B_2$.
\end{corollary}

By making some adjustments in the proof of Theorem \ref{th2}, we can call upon a celebrated theorem by Montgomery and Vaughan \cite{MoVa75} to improve on Corollary \ref{c2}. 

\begin{corollary}\label{c3}
  Suppose that $\alpha_1, \beta_1, \alpha_2, \beta_2$ are real numbers subject to the hypotheses of Theorem \ref{th2}. Then there exists an $\eps = \eps(\bm{\alpha}) > 0$ such that all but $O( x^{1 - \eps} )$ even integers $n \le x$ can be represented as sums of a prime $p_1 \in \mathcal B_1$ and a prime $p_2 \in \mathcal B_2$.
\end{corollary}

Comparing Theorems \ref{th1} and \ref{th2} with the main results in \cite{BaGuNe07}, one notes that our theorems include no hypotheses similar to the condition $\alpha < k$ required in \cite{BaGuNe07}. The latter condition is necessary in the case $\alpha_1 = \dots = \alpha_k = \alpha$, if all large integers $n \equiv k \pmodulo 2$ are to be represented. However, it can be dispensed with when some pair $\alpha_i, \alpha_j$ is linearly independent over $\mathbb Q$.

It seems that the natural hypotheses for the above theorems are that all $\alpha_i$'s and some ratio $\alpha_i/\alpha_j$ be irrational, but such generality is beyond the reach of our method. The finite type conditions above approximate these natural hypotheses without being too restrictive. For example, by a classical theorem of Khinchin's \cite{Khin26}, almost all (in the sense of Lebesgue measure) real numbers are of a finite type.

\section{Preliminaries}
\label{s2}

\subsection*{Notation} 

For a real number $\theta$, $[\theta], \{ \theta \}$ and $\| \theta \|$ denote, respectively, the integer part of $\theta$, the fractional part of $\theta$ and the distance from $\theta$ to the nearest integer; also, $e(\theta) = e^{2\pi i\theta}$. For integers $a$ and $b$, we write $(a, b)$ and $[a, b]$ for the greatest common divisor and the least common multiple of $a$ and $b$. The letter $p$, with or without indices, is reserved for prime numbers. Finally, if $\mathbf x = (x_1, \dots, x_s)$, we write $|\mathbf x| = \max( |x_1|, \dots, |x_s| )$.

\subsection{}\label{s2.1}

For $i = 1, \dots, k$, we set $\gamma_i = \alpha_i^{-1}$ and $\delta_i = \alpha_i^{-1}(1 - \beta_i)$. It is not difficult to see that $m \in \mathcal B_i$ if and only if $0 < \{ \gamma_i m + \delta_i \} < \gamma_i$. Thus, the characteristic function of the Beatty sequence $\mathcal B_i$ is $g_i(\gamma_im + \delta_i)$, where $g_i$ is the 1-periodic extension of the characteristic function of the interval $(0, \gamma_i)$.

Our analysis will require smooth approximations to $g_i$. Suppose that $1 \le i \le k$ and $\Delta$ is a real such that
\[
  0 < \Delta < \sfrac 14 \min (\gamma_i, 1 - \gamma_i).
\]
Then there exist 1-periodic $C^{\infty}$-functions $g_i^{\pm}$ such that:
\begin{itemize}
  \item [i)] $0 \le g_i^-(x) \le g_i(x) \le g_i^+(x) \le 1$ for all real $x$;
  \item [ii)] $g_i^{\pm}(x) = g_i(x)$ when $\Delta \le x \le \gamma_i - \Delta$ or $\gamma_i + \Delta \le x \le 1 - \Delta$.
\end{itemize}
Furthermore, because of ii) above and the smoothness of $g_i^{\pm}$, the Fourier coefficients $\hat g_i^{\pm}(m)$ of $g_i^{\pm}$ satisfy the bounds
\begin{equation}\label{6}
  \hat g_i^{\pm}(0) = \gamma_i + O(\Delta), \qquad |\hat g_i^{\pm}(m)| \ll_r \frac {\Delta^{1 - r}}{(1 + |m|)^r} \quad ( r = 1, 2, \dots).
\end{equation}

\subsection{}\label{s2.2}

The proofs of Theorems \ref{th1} and \ref{th2} use the following generalization of the classical bound for exponential sums over primes in Vaughan \cite[Theorem 3.1]{Vaug97}.

\begin{lemma}\label{l1}
  Suppose that $\alpha$ is real and $a, q$ are integers, with $(a, q) = 1$ and $q \le N$. Then
  \[
    \sum_{p \le N} (\log p)e(\alpha p) \ll \left( Nq^{-1} + N^{4/5} + q \right) \left( 1 + q^2|\theta| \right)(\log 2N)^4,
  \]
  where $\theta = \alpha - a/q$.
\end{lemma}

The proof of the above lemma is essentially the same as that of \cite[Theorem 3.1]{Vaug97}, which is the case $|\theta| \le q^{-2}$. The only adjustment one needs to make in the argument in \cite{Vaug97} is to replace \cite[Lemma 2.2]{Vaug97} by the following variant.

\begin{lemma}\label{l2}
  Suppose that $\alpha, X, Y$ are real with $X \ge 1$, $Y \ge 1$, and $a, q$ are integers with $(a, q) = 1$. Then
  \[
    \sum_{x \le X} \min \left( XYx^{-1}, \| \alpha x \|^{-1} \right) \ll \left( XYq^{-1} + X + q \right) \left( 1 + q^2|\theta| \right)(\log 2Xq) ,
  \]
  where $\theta = \alpha - a/q$.
\end{lemma}

\subsection{}\label{s2.3}

In the next lemma, we use the finite type of an $s$-tuple $\theta_1, \dots, \theta_s$ to obtain rational approximations to linear combinations of $\theta_1, \dots, \theta_s$. 

\begin{lemma}\label{l3}
  Suppose that the $s$-tuple $\theta_1, \dots, \theta_s$ has a finite type and let $\eta > 1$ be such that \eqref{3} has finitely many solutions. Let $0 < \eps < (2\eta)^{-1}$, let $Q$ be sufficiently large, and let $\mathbf m = (m_1, \dots, m_s) \in \mathbb Z^s$, with $0 < |\mathbf m| \le Q^{\eps}$. Then there exist integers $a$ and $q$ such that
  \[
    |q(m_1\theta_1 + \dots + m_s\theta_s) - a| \le Q^{-1}, \quad Q^{\eps} \le q \le Q, \quad (a, q) = 1.
  \]
\end{lemma}

\begin{proof}
  Since the sum $m_1\theta_1 + \dots + m_s\theta_s$ is irrational, it has an infinite continued fraction. Let $q$ and $q'$ be the denominators of two consecutive convergents to that continued fraction, such that $q \le Q < q'$. If $1 \le q \le Q^{\eps}$, then by the properties of continued fractions,
  \begin{align*}
    \|q(m_1\theta_1 + \dots + m_s\theta_s)\| &\le (q')^{-1} < Q^{-1} \le ( q|\mathbf m| )^{-1/(2\eps)}.
  \end{align*}
  When $Q$ is sufficiently large, this contradicts the choice of $\eta$ and $\eps$.
\end{proof}

\section{Proof of Theorem \ref{th1}}
\label{s3}

Without loss of generality, we may assume that the pair of a finite type in the hypotheses of the theorem is $\alpha_1^{-1}, \alpha_2^{-1}$. We also note that if $\alpha_i$ has a finite type, then so does $\gamma_i = \alpha_i^{-1}$. Recall the functions $g_i$ and $g_i^{\pm}$ described in \S\ref{s2.1}. We shall use those functions with $\Delta$ given by
\begin{equation}\label{7}
  \Delta = (\log n)^{-A}.
\end{equation}
We have
\[
  R(n) = \sum_{p_1 + \dots + p_k = n} (\log p_1) \cdots (\log p_k)g_1(\gamma_1p_1 + \delta_1) \cdots g_k(\gamma_kp_k + \delta_k),
\]
so by the construction of $g_i^{\pm}$,
\begin{equation}\label{8}
  R^-(n) \le R(n) \le R^+(n),
\end{equation}
where
\[
  R^{\pm}(n) = \sum_{p_1 + \dots + p_k = n} (\log p_1) \cdots (\log p_k)g_1^{\pm}(\gamma_1p_1 + \delta_1) \cdots g_k^{\pm}(\gamma_kp_k + \delta_k).
\]
We now proceed to evaluate the sums $R^+(n)$ and $R^-(n)$. We shall focus on $R^+(n)$, the evaluation of $R^-(n)$ being similar.

Substituting the Fourier expansions of $g_1^+, \dots, g_k^+$ into the definition of $R^+(n)$, we obtain
\begin{equation}\label{9}
  R^+(n) = \sum_{\mathbf m \in \mathbb Z^k} \hat g_1^+(m_1) \cdots \hat g_k^+(m_k) e(\delta_1m_1 + \dots + \delta_km_k)R(n, \mathbf m),
\end{equation}
where $\mathbf m = (m_1, \dots, m_k)$ and
\[
  R(n, \mathbf m) = \sum_{p_1 + \dots + p_k = n} (\log p_1) \cdots (\log p_k) e( \gamma_1m_1p_1 + \dots + \gamma_km_kp_k ).
\]
We note for the record that when $\mathbf m = \mathbf 0$, we have
\begin{equation}\label{10}
  R(n, \mathbf{0}) = \frac {\mathfrak S_k(n)n^{k - 1}}{(k - 1)!} + O \big( n^{k - 1}(\log n)^{-A} \big).
\end{equation}
When $k = 3$, this is due to Vinogradov \cite{IVin37a} (see also Vaughan \cite[Theorem 3.4]{Vaug97}), and the result for $k \ge 4$ can be proved similarly (see Hua \cite{Hua65}). 

We now set
\begin{equation}\label{11}
  M = \Delta^{-1}(\log n) = (\log n)^{A + 1}.
\end{equation}
Combining \eqref{6}, \eqref{7}, \eqref{9} and \eqref{10}, we deduce that
\begin{equation}\label{12}
  R^+(n) = \frac {\gamma_1 \cdots \gamma_k}{(k - 1)!}\mathfrak S_k(n)n^{k - 1} + O \big( n^{k - 1}(\log n)^{-A} \big) + O( \Sigma_1 + \Sigma_2 ),
\end{equation}
where 
\begin{align*}
  \Sigma_1 &= \sum_{0 < |\mathbf m| \le M} |\hat g_1^+(m_1) \cdots \hat g_k^+(m_k)| |R(n, \mathbf m)|, \\
  \Sigma_2 &= \sum_{|\mathbf m| > M} |\hat g_1^+(m_1) \cdots \hat g_k^+(m_k)| |R(n, \mathbf m)|.
\end{align*}
We may use \eqref{6} to estimate $\Sigma_2$.  It follows easily from the second bound in \eqref{6} that
\begin{equation}\label{13}
  \sum_{m \in \mathbb Z} |\hat g_i^+(m)| \ll \sum_{|m| \le M} \frac 1{1 + |m|} + \sum_{|m| > M} \frac {\Delta^{-1}}{(1 + |m|)^2} \ll \log M.
\end{equation}
Hence, by \eqref{10}, \eqref{11} and \eqref{6} with $r = [A] + 3$,
\begin{align}\label{14}
  \Sigma_2 &\ll R(n, \mathbf 0)(\log M)^{k - 1} \sum_{1 \le i \le k} \sum_{|m| > M} |\hat g_i^+(m)| \\
  &\ll n^{k - 1}(\log n)(\Delta M)^{1 - r} \ll n^{k - 1}(\log n)^{-A}. \notag
\end{align}

Next, we use a variant of the circle method to bound $|R(n, \mathbf m)|$ when $0 < |\mathbf m| \le M$. Define the exponential sum
\[
  S(\xi) = \sum_{p \le n} (\log p)e(\xi p).
\]
By orthogonality,
\begin{equation}\label{15}
  R(n, \mathbf m) = \int_0^1 S(\xi + \gamma_1m_1) \cdots S(\xi + \gamma_km_k) e(-n\xi) \, d\xi.
\end{equation}
Put
\begin{equation}\label{16}
  P = (\log n)^{2A + 12}, \quad Q = nP^{-1}.
\end{equation}
For $j = 1, \dots, k$, we write $\lambda_j = \lambda_j(\mathbf m) = \gamma_jm_j - \gamma_1m_1$. Then
\[
  R(n, \mathbf m) = e(\gamma_1m_1n) \int_{1/Q}^{1 + 1/Q} S(\xi)S(\xi + \lambda_2) \cdots S(\xi + \lambda_k) e(-n\xi) \, d\xi.
\]
We partition the interval $[1/Q, 1 + 1/Q)$ into Farey arcs of order $Q$ and write $\mathfrak M(q, a)$ for the arc containing $a/q$. Thus,
\begin{equation}\label{17}
  |R(n, \mathbf m)| \le \sum_{q \le Q} \sum_{\substack{1 \le a \le q\\ (a, q) = 1}} \int_{\mathfrak M(q, a)} |S(\xi)S(\xi + \lambda_2) \cdots S(\xi + \lambda_k)| \, d\xi.
\end{equation}

When $\xi \in \mathfrak M(q, a)$, with $P < q \le Q$, Lemma \ref{l1} yields
\[
  S(\xi) \ll nP^{-1/2}(\log n)^4 \ll n(\log n)^{-A - 2}.
\]
Inserting this bound into the right side of \eqref{17}, we obtain
\begin{equation}\label{18}
  |R(n, \mathbf m)| \le \sum_{q \le P} \sum_{\substack{1 \le a \le q\\ (a, q) = 1}} \int_{\mathfrak M(q, a)} |S(\xi)S(\xi + \lambda_2) \cdots S(\xi + \lambda_k)| \, d\xi + \Sigma_3,
\end{equation}
where
\begin{align}\label{19}
  \Sigma_3 &\ll n(\log n)^{-A - 2} \sum_{P < q \le Q} \sum_{\substack{1 \le a \le q\\ (a, q) = 1}} \int_{\mathfrak M(q, a)} |S(\xi + \lambda_2) \cdots S(\xi + \lambda_k)| \, d\xi \\
  &\ll n^{k - 2}(\log n)^{-A - 2} \int_0^1 |S(\xi + \lambda_2) S(\xi + \lambda_3)| \, d\xi. \notag
\end{align}
By the Cauchy--Schwarz inequality and Parseval's identity,
\begin{equation}\label{20}
  \int_0^1 |S(\xi + \lambda_i) S(\xi + \lambda_j)| \, d\xi \ll n(\log n) \quad (1 \le i, j \le k), 
\end{equation}
so we deduce from \eqref{19} that
\begin{equation}\label{21}
  \Sigma_3 \ll n^{k - 1}(\log n)^{-A - 1}.
\end{equation}
To estimate the remaining sum on the right side of \eqref{18}, we consider separately the cases $m_1 = 0$ and $m_1 \ne 0$.

\paragraph*{\em Case 1: $m_1 = 0$.}

Since $\mathbf m \ne \mathbf 0$, we have $m_i \ne 0$ for some $i = 2, \dots, k$. By Lemma \ref{l3} with $s = 1$ and $\theta_1 = \gamma_i$, there exist an $\eps > 0$ and integers $b$ and $r$ such that 
\[
  \big| r(m_i\gamma_i) - b \big| < Q^{-1/2}, \quad Q^{\eps} \le r \le Q^{1/2}, \quad (b, r) = 1.
\]
Suppose that $\xi \in \mathfrak M(q, a)$, where $1 \le q \le P$. It follows that
\[
  \left| \xi + \lambda_i - \frac aq - \frac br \right| \le \frac 1{qQ} + \frac 1{rQ^{1/2}} \le \frac 2{rQ^{1/2}}.
\]
Let the integers $a_1, q_1$ be such that
\[
  \frac {a_1}{q_1} = \frac aq + \frac br, \quad (a_1, q_1) = 1.
\]
Then $q_1$ divides $[q, r]$ and is divisible by $[q, r]/(q, r)$, so
\[
  \left| \xi + \lambda_i - \frac {a_1}{q_1} \right| \le \frac {2P}{q_1Q^{1/2}}, \quad rP^{-1} \le q_1 \le rP, \quad (a_1, q_1) = 1.
\]
Thus, Lemma \ref{l1} yields
\begin{align}\label{22}
  S(\xi + \lambda_i) &\ll \big( nq_1^{-1/2} + n^{4/5} + q_1 \big) \left( 1 + q_1PQ^{-1/2} \right) (\log n)^4 \\
  &\ll \big( nq_1^{-1/2} + n^{4/5}P \big) (\log n)^4 \ll n^{1 - \eps/3}. \notag
\end{align}
Combining \eqref{22} and \eqref{20}, we easily get
\begin{equation}\label{23}
  \sum_{q \le P} \sum_{\substack{1 \le a \le q\\ (a, q) = 1}} \int_{\mathfrak M(q, a)} |S(\xi)S(\xi + \lambda_2) \cdots S(\xi + \lambda_k)| \, d\xi \ll n^{k - 1 - \eps/4}.
\end{equation}

\paragraph*{\em Case 2: $m_1 \ne 0$.}

Then we apply Lemma \ref{l3} to the pair $\gamma_1, \gamma_2$. It follows that there exist an $\eps > 0$ and integers $b$ and $r$ such that
\[
  \big| r(m_2\gamma_2 - m_1\gamma_1) - b \big| < Q^{-1/2}, \quad Q^{\eps} \le r \le Q^{1/2}, \quad (b, r) = 1.
\]
Suppose that $\xi \in \mathfrak M(q, a)$, where $1 \le q \le P$. Arguing similarly to Case 1, we find that there exist integers $a_1, q_1$ such that
\[
  \left| \xi + \lambda_2 - \frac {a_1}{q_1} \right| \le \frac {2P}{q_1Q^{1/2}}, \quad rP^{-1} \le q_1 \le rP, \quad (a_1, q_1) = 1.
\]
Using this rational approximation to $\xi + \lambda_2$, we can now apply Lemma \ref{l1} to show that
\[
  S(\xi + \lambda_2) \ll \big( nq_1^{-1/2} + n^{4/5}P \big) (\log n)^4 \ll n^{1 - \eps/3}. 
\]
We then derive \eqref{23} in a similar fashion to Case 1.

We conclude that \eqref{23} holds for all vectors $\mathbf m$ with $0 < |\mathbf m| \le M$. Together, \eqref{18}, \eqref{21} and \eqref{23} yield
\[
  |R(n, \mathbf m)| \ll n^{k - 1}(\log n)^{-A - 1} 
\]
for all $0 < |\mathbf m| \le M$, whence
\begin{equation}\label{24}
  \Sigma_1 \ll n^{k - 1}(\log n)^{-A - 1}(\log M)^k \ll n^{k - 1}(\log n)^{-A}.
\end{equation}
Finally, from \eqref{12}, \eqref{14} and \eqref{24},
\[
  R^+(n) = \frac {\gamma_1 \cdots \gamma_k}{(k - 1)!}\mathfrak S_k(n)n^{k - 1} + O \big( n^{k - 1}(\log n)^{-A} \big).
\]
Since an analogous asymptotic formula holds for $R^-(n)$, the conclusion of the theorem follows from \eqref{8}. \qed

\section{Sketch of the proof of Theorem \ref{th2}}
\label{s4}

Let $R^+(n)$ and $R^-(n)$ be the quantities defined in \S\ref{s3} with $k = 2$. To prove Theorem \ref{th2} it suffices to establish the inequality
\begin{equation}\label{25}
  \sum_{n \le x} \left| R^{\pm}(n) - \gamma_1\gamma_2\mathfrak S_2(n)n \right|^2 \ll x^3(\log x)^{-3A}.
\end{equation}
As in the proof of Theorem \ref{th1}, we focus on the proof of the inequality for $R^+(n)$, the proof of the other inequality being similar.

We use the notation introduced in \S\ref{s3} with $k = 2$ and $A$ replaced by $2A + 1$. When $k = 2$, the asymptotic formula \eqref{10} is not known, but we do have the upper bound (see \cite{HaRi74})
\[
  R(n, \mathbf 0) \ll n\prod_{p \mid n} \left( \frac p{p - 1} \right) \ll n\log\log n.
\] 
This bound suffices to show similarly to \eqref{9}--\eqref{14} that 
\[
  R^+(n) = \gamma_1\gamma_2R(n, \mathbf 0) + O\left( n(\log n)^{-2A} + \Sigma_1 \right).
\]
Furthemore, by \cite[Theorem 3.7]{Vaug97},
\[
  \sum_{n \le x} \left| R(n, \mathbf 0) - \mathfrak S_2(n)n \right|^2 \ll x^3(\log x)^{-3A}.
\]
Thus, \eqref{25} for $R^+(n)$ follows from the inequality
\begin{equation}\label{26}
  \sum_{n \le x} \bigg| \sum_{0 < |\mathbf m| \le M} |\hat g_1^+(m_1) \hat g_2^+(m_2)| |R(n, \mathbf m)| \bigg|^2 \ll x^3(\log x)^{-3A}.
\end{equation}
By \eqref{13} and Cauchy's inequality, the left side of \eqref{26} is
\[
  \ll (\log M)^4 \max_{0 < |\mathbf m| \le M} \sum_{n \le x} |R(n, \mathbf m)|^2,
\]
so it suffices to show that
\begin{equation}\label{27}
  \sum_{n \le x} |R(n, \mathbf m)|^2 \ll x(\log x)^{-3A - 1}
\end{equation}
for all $\mathbf m$, $0 < |\mathbf m| \le M$. By \eqref{15} and Bessel's inequality,
\begin{equation}\label{28}
  \sum_{n \le x} |R(n, \mathbf m)|^2 \le \int_0^1 |S(\xi + \gamma_1m_1)S(\xi + \gamma_2m_2)|^2 \, d\xi,
\end{equation}
where the definition of the exponential sum $S(\xi)$ has been altered to
\[
  S(\xi) = \sum_{p \le x} (\log p)e(\xi p).
\]
We set
\begin{equation}\label{29}
  P = (\log x)^{3A + 10}, \quad Q = xP^{-1},
\end{equation}
and obtain similarly to \eqref{17} that
\begin{equation}\label{30}
  \int_0^1 |S(\xi + \gamma_1m_1)S(\xi + \gamma_2m_2)|^2 \, d\xi \le \sum_{q \le Q} \sum_{ \substack{ 1 \le a \le q\\ (a, q) = 1}}  \int_{\mathfrak M(q, a)} |S(\xi)S(\xi + \lambda_2)|^2 \, d\xi.
\end{equation}
As in \S\ref{s3},
\[
  \min \left( |S(\xi)|, |S(\xi + \lambda_2)| \right) \ll xP^{-1/2}(\log x)^4
\]
for all $\xi$, so we deduce from \eqref{20}, \eqref{29} and \eqref{30} that
\[
  \int_0^1 |S(\xi + \gamma_1m_1)S(\xi + \gamma_2m_2)|^2 \, d\xi \ll x^3P^{-1}(\log x)^9 \ll x^2(\log x)^{-3A - 1}.
\]
Inserting the last bound into \eqref{28}, we obtain \eqref{27}. \qed

\section{Closing remarks}
\label{s5}

In the proof of Theorem \ref{th1}, we essentially showed that
\begin{equation}\label{31}
  R(n) = \gamma_1 \cdots \gamma_k R(n, \mathbf 0) + \text{error terms},
\end{equation}
and then chose the parameters $\Delta, M, P, Q$ so that the error terms were $\ll n^{k - 1}(\log n)^{-A}$. It is possible to alter the above choices so that the error terms in \eqref{31} are $\ll n^{k - 1 - \eps}$ for some $\eps > 0$ which depends only on the $\alpha_i$'s. Therefore, the quality of the error term in \eqref{4} is determined solely by the quality of the error term in the asymptotic formula \eqref{10} for the number of representations of an integer $n$ as the sum of $k$ primes. However, since no improvements on \eqref{10} are known, the improved bounds for the error terms in \eqref{31} have no effect on Theorem \ref{th1}. 

Similarly, when $k = 2$, a slight alteration of our choices in \S\ref{s4} yields the bound
\begin{equation}\label{32}
  R(n) = \gamma_1\gamma_2 R(n, \mathbf 0) + O\big( n^{1 - \eps} \big)
\end{equation}
for all but $O( x^{1 - \eps} )$ values of $n \le x$. In this case, however, such a variation has a tangible effect: it yields Corollary \ref{c3}. Indeed, by a well-known result of Montgomery and Vaughan~\cite{MoVa75}, there is an absolute constant $\omega < 1$ such that the right side of \eqref{32} is positive for all but $O(x^{\omega})$ even integers $n \le x$. 

Finally, a comment regarding our finite type hypotheses. We say that an $s$-tuple $\theta_1, \dots, \theta_s$ of real numbers \emph{is of subexponential type}, if for each fixed $\eta > 0$, the inequality
\[
  \| q_1\theta_1 + \dots + q_s\theta_s \| < \exp ( -|\mathbf q|^{\eta} )
\]
has only finitely many solutions $\mathbf q = (q_1, \dots, q_s) \in \mathbb Z^s$. Clearly, every $s$-tuple of a finite type is also of subexponential type, but not vice versa. It takes little effort to check that in the arguments in \S\ref{s3} and \S\ref{s4}, it suffices to assume that each $\alpha_i$ and some pair $\alpha_i, \alpha_j$ are of subexponential type. Thus, our method reaches some Beatty sequences $\mathcal B_{\alpha, \beta}$ with $\alpha$ of an infinite type. On the other hand, under the weaker subexponential type hypotheses, we no longer have the improved remainder estimates in \eqref{31} and \eqref{32}. In particular, we no longer have Corollary \ref{c3} (at least, not by the simple argument sketched above). This seems to be too steep a price to pay for such a modest gain in generality.

\bibliographystyle{amsplain}

\end{document}